\documentclass[10pt]{amsart}
\usepackage{amssymb}

\usepackage{amscd}
\usepackage{amsthm}

\usepackage[all]{xy}
\usepackage[dvips]{graphicx}

\newtheorem{Thm}{Theorem}[section]
\newtheorem{Lem}[Thm]{Lemma}

\newtheorem{Prop}[Thm]{Proposition}
\newtheorem{Ex}[Thm]{Example}
\newtheorem{Rem}[Thm]{Remark}
\newtheorem{Conj}[Thm]{Conjecture}

\newenvironment{Proof}{\noindent{\it Proof.\ \
}}{\phantom{}\hfill$\square$}

\begin{document}

\title[Vertices of indecomposable modules]
{On the vertices of indecomposable modules over   dihedral
$2$-groups}

\author{Guodong ZHOU}
\address{Guodong ZHOU\newline
 LAMFA et CNRS UMR 6140 \newline Universit\'e de Picardie Jules Verne\newline
33, rue St Leu\newline 80039 Amiens\newline FRANCE}
\email{guodong.zhou@u-picardie.fr}


\begin{abstract}Let $k$ be an algebraically closed field of
characteristic $2$.  We calculate the vertices of all indecomposable
$kD_8$-modules  for   the dihedral group $D_8$ of order $8$. We¡¡
also give a conjectural formula of the induced module
  of a string module   from $kT_0$ to $kG$
where $G$ is    a dihedral group $G$ of order $\geq 8$  and where
$T_0$ is a dihedral subgroup of index $2$ of $G$. Some cases where
we verified this formula are given.

\end{abstract}

\maketitle





\section{Introduction}



 Let $k$ be an algebraically closed
 field  and let $G$ be a finite group. A subgroup $D\leqslant G$ is called a  vertex   of an indecomposable $kG$-module
 $M$¡¡  if
    $M$ is a direct summand of an induced module from $D$ to
 to   $G$ and if $D$ is minimal for this property. It can be easily seen that two vertices of
 $M$ are conjugate in $G$. We  denote by $vx(M)$ a vertex of
 $M$. The knowledge of vertices of modules is a
central point in modular representation theory of finite groups In
particular, it is important   to understand the
 category of modules over   group algebras. To determine the vertex
of an indecomposable module   usually is  a hard problem. Much work
has been done on this problem and  is mainly  centered around  the
vertex of a simple module (see \cite{HM76},\cite{K79} for general
statements).
 According to a theorem of V.M.Bondarenko and Yu.A.Drozd(\cite{BD77}),  a block $B$ of a group
 algebra has finitely many isomorphism classes of indecomposable
 modules (i.e. $B$  has finite representation type) if and only if it has   cyclic defect groups. So
blocks with cyclic defect groups  are the  easiest to study,  see
for example \cite{D66}\cite{J69}\cite{P75}\cite{M76}\cite{M75}. The
case of tame representation type is a natural continuation to deal
with and here by the classification of tame blocks those of dihedral
defect groups are natural candidates. Only special situations are
known. We mention a few of them. Vertices of simple modules for the
case of blocks with cyclic defect groups were calculated in
\cite{M76} ( and also vertices of all indecomposable modules in
\cite{M75}). K. Erdmann dealt with some blocks with dihedral defect
groups in \cite{E77}. Some other group algebras of not necessarily
tame representation type are also considered in the literature, see
\cite{LM78} \cite{Wil}\cite{MZ} etc.


  In this note, we treat the dihedral group of
order eight. Let $k$ be an algebraically closed field of
characteristic $2$. Given  $D_8$ the
 dihedral defect group of order $ 8$,
 using purely linear
 algebra method,   we
compute  the induced modules of all indecomposable modules from each
subgroup to $D_8$ and we thus  obtain the   vertices of \textit{all}
indecomposable modules. Roughly speaking, in the Auslander-Reiten
quiver of $kD_8$, for a homogenous tube, all modules have the same
vertex, or the module at the bottom has a smaller vertex and all
other modules have the same; for a component of
$\mathbb{Z}A^{\infty}_{\infty}$ type, if different vertices appear,
there are two $\tau$-orbits which have the same vertex and all the
other modules have a larger vertex.

Since the pioneering work of P.Webb (\cite{W82}), the relations
between inductions from subgroups and the Auslander-Reiten quiver
are extensively studied, see \cite{E86}\cite{K89}\cite{K1}\cite{K2}.
The distribution of vertices of modules in the Auslander-Reiten
quiver becomes an interesting problem. This problem was solved in
case of  $p$-groups in \cite{E90} and in \cite{U91}\cite{OU94} in
  general.  In fact, K.Erdmann considered all components except for homogenous tubes in the Auslander-Reiten quiver.
 Our results  thus verify  her result   in the special case of the dihedral group of order
$8$ and furthermore complement it by
  dealing with homogenous tubes which cause most of the  difficulties of calculations.

  For dihedral $2$-groups of order $\geq 16$, we only obtain partial
  results, but we propose a conjectural formula for the induced
  module of a string module from a dihedral subgroup of index $2$ to
  the whole dihedral group. This formula should give the vertices
  of all string modules. More precisely,
let $$G=D_{2^n}=\langle x, y|x^{2^{n-2}}=e=y^2, yxy=x^{-1}\rangle$$
be the dihedral group of order $2^n$ with $n\geq 3$ and let
$T_0=\langle x^2, y\rangle$ be a dihedral group of index $2$. Let
$M(C)$ be a string module over $kT_0$ (for the definition of a
string module, see Section ~\ref{mr}). Then we construct a new
string $\varphi(C)$  over $kG$ (for details see Section 4) and  the
following formula should hold
$$\text{Ind}_{T_0}^GM(C)\cong M(\varphi(C)).$$

 This paper is organized as follows. In Section \ref{mr} we
present the classification of indecomposable modules over dihedral
$2$-groups. Vertices of indecomposable modules over the dihedral
 group of order eight are calculated in the third section, where the   main theorem of this
paper Theorem~\ref{main} is prove, but we postpone  in  the final
section the proof of Proposition \ref{longcorde} which is rather
technical in nature.
   We give the induction formula in Section~\ref{formula} and some
special cases of this formula are proved.

\textbf{Notations}: $\mathbb{N}_0=\{0, 1, 2, \cdots\}$ and
$\mathbb{N}_1=\{1, 2, 3,  \cdots\}$

\textbf{Acknowledgement}: This paper is part of my Ph. D thesis
defended on the 19 June 2007. I want to express my sincere gratitude
to my thesis supervisor Prof. Alexander Zimmermann for his patient
guidance and his constant encouragements.


 \section{Classification of indecomposable modules over dihedral $2$-groups}\label{mr}

Since the pioneering work of P. Gabriel (\cite{G}), quivers become
important in representation theory. A theorem of P. Gabriel says
that any finite-dimensional algebra over an algebraically closed
field is Morita equivalent to an algebra, its basic algebra, defined
by quiver with relations. We will use  the presentation by quiver
with relations throughout the present note. For the general theory
of quiver with relations, see \cite[Chapter 4]{Ben} or \cite[Chapter
3]{ARS}.

 Let
$$G=D_{2^n}=\langle x, y|x^{2^{n-1}}=e=y^2, yxy=x^{-1}\rangle$$ be the dihedral
group of order $2^n$ with $n\geq 2$. The group algebra $kG$ is basic
  and its quiver with relations  is of the
following form:
 $$\unitlength0.6cm
\begin{picture}(7,3)
 \put(5.4,2){\circle{2.0}}
\put(2.8,2){\circle{2.0}} \put(3.85,1.9){\vector(0,1){0.3}}
\put(3.95,1.9){$\bullet$}  \put(4.35,1.9){\vector(0,1){0.3}}
\put(0.8,1.9){$\alpha$} \put(7,1.9){$\beta$}
 \put(0.8,0){$ \alpha^2=0=\beta^2, \ \ (\alpha \beta)^{2^{n-2}}=(\beta \alpha)^{2^{n-2}}$}

\end{picture}$$
where $\alpha=1+y$ and $\beta=1+yx$.

For the convenience of later use, we record some subgroups of $G$
and the quiver with relations of the corresponding group algebras.
Note that
$$H=\langle
x\rangle,\ \ T_0=\langle x^2, y\rangle, \ \ T_1=\langle x^2,
yx\rangle$$
 These
are all the subgroups  of index $2$ of $G$. Furthermore,  $H \simeq
C_{2^{n-2}}$ is the cyclic group of order $2^{n-2}$ et $T_0\simeq
D_{2^{n-1}} \simeq T_1$ are isomorphic to  the dihedral group of
order $2^{n-1}$. The quiver with relations of $kT_0$, $kT_1$ and
$kH$ are respectively
$$\unitlength0.6cm
\begin{picture}(7,3)
 \put(5.4,2){\circle{2.0}}
\put(2.8,2){\circle{2.0}} \put(3.85,1.9){\vector(0,1){0.3}}
\put(3.95,1.9){$\bullet$}  \put(4.35,1.9){\vector(0,1){0.3}}
\put(0.8,1.9){$\alpha_0$} \put(7,1.9){$\beta_0$}
 \put(0.8,0){$ \alpha_0^2=0=\beta_0^2,\  (\alpha_0 \beta_0)^{2^{n-3}}= (\beta_0 \alpha_0)^{2^{n-3}}$}

\end{picture}$$
where $\alpha_0=1+y$ and $\beta_0=1+yx^2$,
$$\unitlength0.6cm
\begin{picture}(7,3)
 \put(5.4,2){\circle{2.0}}
\put(2.8,2){\circle{2.0}} \put(3.85,1.9){\vector(0,1){0.3}}
\put(3.95,1.9){$\bullet$}  \put(4.35,1.9){\vector(0,1){0.3}}
\put(0.8,1.9){$\alpha_1$} \put(7,1.9){$\beta_1$}
 \put(0.8,0){$ \alpha_1^2=0=\beta_1^2,\  (\alpha_1 \beta_1)^{2^{n-3}}= (\beta_1
 \alpha_1)^{2^{n-3}}$}

\end{picture}$$
where $\alpha_1=1+yx$ and $\beta_1=1+ yx^3$ and
$$\unitlength0.6cm
\begin{picture}(7,2)
 \put(2.4,1){\circle{2.0}}
\put(1,0.9){$\bullet$} \put(3.8,0.9){$\gamma$}

\put(1.35,0.9){\vector(0,1){0.3}} \put(5,0.9){$\gamma^{2^{n-1}}=0$}
\end{picture}$$
where $\gamma=1+x $.

 Inspired by the work  \cite{GP} of I.M.Gelfand and
V.A.Ponomarev on the representation theory of the Lorentz group,
C.M.Ringel in \cite{Rin} classified  indecomposable modules over
$kG$. Precisely except the module of the entire group algebra $kG$,
all other indecomposable modules can be divided into two families:
string modules and band modules. We now recall his classification.

We define two strings $1_{\alpha}$ and $1_{\beta}$ of length zero
with $1_{\alpha}^{-1}=1_{\beta}$ and  $1_{\beta}^{-1}=1_{\alpha}$.
Consider now $\alpha, \beta, \alpha^{-1}, \beta^{-1}$ as 'letters'
in formal language and let $(\alpha^{-1})^{-1}=\alpha$ and
$(\beta^{-1})^{-1}=\beta$. If $l$ is a letter, we write $l^*$ to
mean   'either $l$ or $l^{-1}$'. A \textit{string} $C=l_1l_2\cdots
l_n$ of length $n\geq 1$ is given by a sequence $l_1l_2\cdots l_n$
of letters subject to
\begin{itemize}
\item $l_i=\alpha^*$ for $1\leq i\leq n-1$ implies $l_{i+1}=\beta^*$
and similarly $l_i=\beta^*$ for $1\leq i\leq n-1$ implies
$l_{i+1}=\alpha^*$

\item for any $1\leq i<j\leq n$, neither $l_i\cdots l_j$ nor $l_j^{-1}\cdots l_i^{-1}$ is  in the set
$\{(\alpha\beta)^{2^{n-2}},  (\beta\alpha )^{2^{n-2}}\}$.
\end{itemize}
For instance, the word $C=\alpha\beta^{-1}\alpha^{-1}\beta$ is a
string of length $4$. We illustrate usually this string by the
following graph:
$$\vcenter
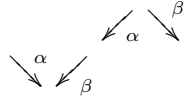
{\xymatrix@R0.4cm@C0.4cm{&&&\ar[dr]^\beta\ar[dl]^\alpha&\\
                     \ar[dr]^\alpha& & \ar[dl]^\beta&&\\
                     &&&&}}$$
In this graph, we draw an arrow from north-west to south-east for a
direct letter, and an arrow  from north-east to south-west for an
inverse letter. If $C=l_1\cdots l_n$ is a string, then its inverse
is given by $C^{-1}=l_n^{-1}\cdots l_1^{-1}$. Let $\mathcal{S}t$ be
the set of all strings. Let $\rho$ be the equivalence relation on
$\mathcal{S}t$ which identifies each string to its inverse. If
$C=l_1\cdots l_n$ and $D=f_1\cdots f_m$ are two strings, their
product is given by $CD=l_1\cdots l_nf_1\cdots f_m$ provided that
this is again a string.
 Let $\mathcal{B}d$ be the set of strings of
even length $\neq 0$ and which are not powers of strings of strictly
smaller length. The elements of $\mathcal{B}d$ are called
\textit{bands}. If $C=l_1\cdots l_n$ is a band, then  for $1\leq
i\leq n-1$ denote by $C_{(i)}$ the $i$-th cyclic permutation word,
thus $C_{(0)}=l_1\cdots l_n$, $C_{(1)}=l_2\cdots l_nl_1$, up to
$C_{(n-1)}=l_nl_1\cdots l_{n-1}$. Let $\rho'$ be the equivalence
relation which identifies with the band $C$ all its  cyclic
permutations $C_{(i)}$ and their inverses $C_{(i)}^{-1}$.

To every string $C$, we are going to construct an  indecomposable
module, denoted by $M(C)$ and called a \textit{string module}.
Namely, let $C=l_1\cdots l_n$ be a string of length $n$. Let $M(C)$
be given by a $K$-vector space of dimension $n+1$, say with basis
$e_0, e_1, \cdots, e_n$ on which $\alpha$ and $\beta$ operate
according to the following schema
$$\xymatrix@1{ke_0 \ar@{-}[r]^{l_1}&  ke_1 \ar@{-}[r]^{l_2} &ke_2 \ar@{-}[r]^{l_3}&
\ar@{.}[r]
  & \ar@{-}[r]^{l_{n-1}}&  ke_{n-1}\ar@{-}[r]^{l_n} &k e_n}.$$
For example, if $C= \alpha\beta^{-1}\alpha^{-1}\beta$, we have the
following schema
$$\vcenter
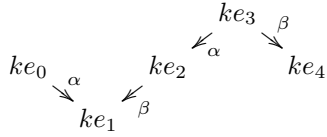
{\xymatrix@R0.2cm@C0.2cm{&&&ke_3\ar[dr]^\beta\ar[dl]^\alpha&\\
                     ke_0\ar[dr]^\alpha& & ke_2\ar[dl]^\beta&&ke_4\\
                     &ke_1&&&}}$$
Note that we already use the notation above to adjust the direction
of the arrows according to whether the letter $l_i$ is direct or
not. This graph indicates how the basis vectors $e_i$ are mapped
into each other or into zero, more precisely,
$$\alpha e_0=e_1, \alpha e_1=0, \alpha e_2=0, \alpha e_3=e_2, \alpha
e_4=0$$ and
$$\beta e_0=0, \beta e_1=0, \beta e_2=e_1, \beta e_3=e_4, \beta
e_4=0.$$ It is obvious that $M(C)$ and $M(C^{-1})$ are isomorphic.

Next we construct   \textit{band modules}. Let $\lambda \in k^*$ and
$n\in \mathbb{N}_1$. For each equivalence class of bands with
respect to $\rho'$, we choose one representative $C=l_1\cdots l_n$
such that $l_n$ is inverse.  Let $M(C, n, \lambda)$ be given by
$M(C, n, \lambda)=\oplus_{i=0}^{n-1}V_i$ with $V_i=k^n$  for any
$1\leq i\leq n$ on which $\alpha$ and $\beta$ operate according to
the following schema
$$
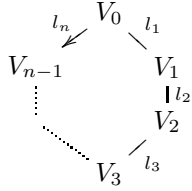
\xymatrix@R0.2cm@C0.2cm{  &  V_0 \ar@{-}[dr]^{l_1} \ar[dl]_-{l_n} & \\
             V_{n-1} \ar@{.}[d] & & V_1 \ar@{-}[d]^{l_2} \\
            \ar@{.}[dr] & & V_2\ar@{-}[dl]^{l_3}\\
            &  V_3&} $$This means that the action is given by
            \begin{itemize}
\item[(1)]$l_s:V_{s-1} \rightarrow V_s$ is the identity map,  if  $l_s$ is direct for   $1
\leq s \leq n-1$,

\item[(2)]$l_s^{-1}:V_s \rightarrow V_{s-1}$ is the identity map,  if $l_s$ is
inverse for $1 \leq s \leq n-1$,

\item[(3)]$l_n^{-1}:V_n=V_0 \rightarrow V_{n-1}$ is $J_n(\lambda)$ where
$J_n(\lambda)$ is the block of Jordan.
\end{itemize}

\begin{Thm}(\cite[Section 8]{Rin})
The strings modules $M(C)$ with $C\in \mathcal{S}t/\rho$ and the
band  modules $M(C, n, \lambda)$ with $C\in \mathcal{B}d/\rho'$,
$n\in \mathbb{N}_1$ and $\lambda\in k^*$,  together with $kG$,
furnish a complete list of isomorphism classes of indecomposable
$kG$-modules for $G$ the dihedral group of order $2^n$ with $n\geq
2$.
\end{Thm}

Now we can  describe the Auslander-Reiten quiver of $kG$. For the
general theory of Auslander-Reiten quiver, we refer to \cite[Chapter
4]{Ben} and \cite{ARS}. Let $C$ be a string, we denote  by $Q(C)$
the component of the Auslander-Reiten quiver containing $M(C)$. Let
$D$ be a band, $n\geq 1$ and $\lambda\in k^*$. Then we denote   by
$Q(D,\lambda)$ the component of the Auslander-Reiten quiver
containing $M(D, n, \lambda)$.

\begin{Prop}\label{klein}(\cite[Chaptre 4, Section 4.17]{Ben}) Let $G$ be the Klein-four group.
 The Auslander-Reiten quiver of $kG$ consists of

\begin{itemize}

\item[$\bullet$] infinitely many  homogenous tubes $Q(\alpha\beta^{-1}, \lambda)$ with $\lambda \in k^*$
formed by band modules $M(\alpha\beta^{-1}, n, \lambda)$ with $n\in
\mathbb{N}_1$,

\item[$\bullet$] two homogeneous tubes
$Q(\alpha ) $ and $Q(\beta ) $,

\item[$\bullet$] one component of type
$\mathbb{Z}\tilde{A}_{12}$   consisting of all the syzygies of the
trivial module of dimension $1$.

\end{itemize}
\end{Prop}
\begin{Prop}\label{dihedral}(\cite[Chaptre 4, Section 4.17]{Ben}) Let $G$ be a dihedral group of order $\geq 8$.
 The Auslander-Reiten quiver of $kG$ consists of

\begin{itemize}

\item[$\bullet$] infinitely many  homogenous tubes
consisting of band modules,

\item[$\bullet$] two homogeneous tubes
$Q((\alpha\beta)^{2^{n-2}-1}\alpha) $ and
$Q((\beta\alpha)^{2^{n-2}-1}\beta) $,

\item[$\bullet$] infinitely many  components of
$\mathbb{Z}A^{\infty}_{\infty}$ type consisting of string modules.

\end{itemize}
\end{Prop}

\section{The dihedral group of order eight}\label{eight}
\subsection{Statement of the main theorem}
Now we specialize to the dihedral group of order $8$. We fix    some
notations.
  Let
$$D_8=\langle x, y|x^4=e=y^2, yxy=x^{-1}\rangle$$ be the dihedral
group of order $8$. The   quiver with relations of $kD_8$ is given
in Section ~\ref{mr}. Note that
$$H=\langle
x\rangle=\{e,x,x^2,x^3\},\ \ T_0=\langle x^2, y\rangle=\{e, y, x^2,
yx^2\}, \ \ T_1=\langle x^2, yx\rangle=\{e, yx, x^2, yx^3\}$$
 Recall that $H
\simeq C_4$ is the cyclic group of order $4$ et $T_0\simeq V_4
\simeq T_1$ are  the Klein-four group. Their quiver with relations
are also given in Section~\ref{mr}.

In order to state the main theorem, we introduce some particular
bands. For $n\leq 2$, we define
$$C_1=\beta\alpha^{-1}\beta^{-1}\alpha=
\vcenter{\xymatrix@R0.2cm@C0.2cm{&&& \ar[dr]^\alpha \ar[dl]^\beta&\\
\ar[dr]^\beta& &\ar[dl]^\alpha&&\\&
 && & & }}$$
$$C_2=\beta\alpha\beta\alpha^{-1}\beta^{-1}\alpha\beta^{-1}\alpha^{-1}=
\vcenter{\xymatrix@R0.2cm@C0.2cm{\ar[dr]^\beta&&&& & & &&
\ar[dl]^\alpha\\
& \ar[dr]^\alpha & & & & \ar[dr]^\alpha \ar[dl]^\beta&& \ar[dl]^\beta& \\
& &  \ar[dr]^\beta& &\ar[dl]^\alpha&&& & & \\&
 && & & & & & & }}$$
If  $n\geq 2$ is even,
 $C_{n+1}=\alpha\beta\alpha^{-1}C_n\beta^{-1}$; if   $n\geq 3$ is odd, $C_{n+1}=\beta
 C_n\alpha\beta^{-1}\alpha^{-1}$.  For $n\in \mathbb{N}_1$, we define $D_n$
 as the band obtained by exchanging
   $\alpha$ and $\beta$ in $C_n$. Notice that $C_1=D_1^{-1}$.

\begin{Thm}\label{main}
 Let  $M$ be an    indecomposable module over the group $kD_8$ where
$k$ is   an algebraically closed
 field of characteristic $2$ and where
  $$D_8=\langle x, y |
x^4=e=y^2, yxy=x^{-1}\rangle$$ is  the dihedral group of order $8$.
Denote by $ vx(M) $ a vertex of $M$. Then
\begin{itemize}
\item[(1)] $vx(kD_8)=\{e\}$.

\item[(2)] In the homogeneous tube $Q(\beta \alpha^{-1}, 1)$,
we obtain $vx(M(\beta \alpha^{-1}, 1, 1))=H=\langle x \rangle$ and
$vx(M)=D_8$ for $M\ncong M(\beta \alpha^{-1}, 1, 1)$.

\item[(3)] In the homogeneous tube  $Q(\beta\alpha\beta^{-1}\alpha^{-1}, 1)$,
we have  $$vx(M(\beta\alpha\beta^{-1}\alpha^{-1}, 1, 1))=\langle
x^2\rangle=\{e,x^2\}$$ and for $M\ncong
M(\beta\alpha\beta^{-1}\alpha^{-1}, 1, 1)$,  $vx(M)=D_8$.

\item[(4)] In the homogeneous tube $Q(\beta\alpha\beta \alpha^{-1}\beta^{-1}\alpha^{-1}, 1)$,
  $vx(M(\beta\alpha\beta \alpha^{-1}\beta^{-1}\alpha^{-1},
1, 1))=H$ and  $vx(M)=D_8$ for $M\ncong M(\beta\alpha\beta
\alpha^{-1}\beta^{-1}\alpha^{-1}, 1, 1)$.

\item[(5)] Each module in the homogeneous tube  $Q(\beta\alpha\beta \alpha^{-1}, \mu)$
 with  $\mu \in k^*$
has the vertex $$T_0=\langle x^2, y \rangle= \{e, y, x^2, yx^2\}.$$

\item[(5')] Each  module in the homogeneous tube $Q(\alpha\beta\alpha\beta^{-1}, \mu)$ with $\mu \in k^*$
has the vertex $$T_1=\langle x^2, yx\rangle=\{e, yx, x^2, yx^3\}.$$

\item[(6)] In the homogeneous tube $Q(C_1, 1)$,
we have $$vx(M(C_1, 1, 1))=\langle x^2\rangle=\{e, x^2\}$$ and for
$M\ncong M(C_1, 1, 1)$,    $vx(M)=D_8$.

\item[(7)] In the homogeneous tube $Q(C_n, 1)$ with $n\geq 2$, we
have  $vx(M(C_n, 1, 1))=T_0$ and for $M\ncong M(C_n, 1, 1)$,
$vx(M)=D_8$.

\item[(7')] In the homogeneous tube $Q(D_n, 1)$ with $n\geq 2$, we have
 $vx(M(D_n, 1, 1))=T_1$ and for $M\ncong M(D_n, 1, 1)$,
$vx(M)=D_8$.

\item[(8)] In the homogeneous tube $Q(\beta\alpha\beta)$,
$vx(M(\beta\alpha\beta))=\langle y\rangle=\{e,y\}$ and for $M\ncong
M(\beta\alpha\beta )$, $vx(M)=T_0$.

\item[(8')] In the homogeneous tube $Q(\alpha\beta\alpha)$,
$vx(M(\alpha\beta\alpha))=\langle yx\rangle=\{e,yx\}$ and for
$M\ncong M(\alpha\beta\alpha )$, $vx(M)=T_1$.

\item[(9)] In the  component $Q( \beta )$ of type $\mathbb{Z}A^{\infty}_{\infty}$, the syzygies $\Omega^n
(M(\beta))$ with $n \in \mathbb{Z}$ which form two $\tau$-orbits
have vertices $T_0$ and all others  modules have $D_8$ as vertices.

\item[(9')] In the  component $Q( \alpha )$ of type $\mathbb{Z}A^{\infty}_{\infty}$, the syzygies  $\Omega^n
(M(\alpha))$ with $n \in \mathbb{Z}$  which form two $\tau$-orbits
have vertices $T_1$ and all others modules have $D_8$ as vertices.

\item[(10)] All other indecomposable modules have    $D_8$ as
vertices.

\end{itemize}
\end{Thm}

We  next  collect some well-known results about cyclic groups and
Klein-four group which will be needed in the proof of the preceding
theorem.

\begin{Lem}\label{cyclic}
Let $H=\langle x\rangle $ be the cyclic subgroup of order $4$ of
$D_8$. Then
\begin{itemize}
\item[(1)] each indecomposable $kH$-module is of the form $M(\gamma^i)$
for $0\leq i\leq 3$.

\item[(2)] we have $vx(M(\gamma^0))=H=vx(M(\gamma^2))$, $vx(M(\gamma^1))=\langle
x^2\rangle$ and $vx(M(\gamma^3))=vx(kH)=\{e\}$.

\end{itemize}
\end{Lem}
\begin{Proof} For (1), see \cite[Section II.4]{A}, we just translate the description there into the context of string modules.
 For (2), it is sufficient to calculate the induced module of the trivial
module from the subgroup $\{e, x^2\}$ to $H$.
\end{Proof}

\begin{Lem}\label{klein1}

Let $T_0=\langle x^2, y\rangle \cong V_4$. Then
\begin{itemize}

\item[(i)] $vx(kT_0)=\{e\}$

\item[(ii)] $vx(M(\beta_0))=\langle y\rangle= \{e, y\}$, $vx(M(\alpha_0))=\langle yx^2\rangle=\{e, yx^2\} $ et $
vx(M(\alpha\beta^{-1},1,1))=\langle x^2\rangle=\{e, x^2\}$

\item[(iii)] any other indecomposable module has    $T_0$ as a vertex.

\end{itemize}

 \end{Lem}
\begin{Proof}It is sufficent to compute inductions from its subgroups to $T_0$.
Recall the general method to calculate  an induced module. Let $G$
be a finite group and $H$ a subgroup of index $m$. Then we write
 $G=\coprod_{i=1}^m g_iH$. For a  $kH$-module $M$, its induced
 module is
$$\text{Ind}_H^GM=\coprod_{i=1}^mg_i \otimes M$$  and  the action
is given by   $g(g_i \otimes x)=g_j\otimes hx, $   for all $g \in
G$,   $1 \leq i \leq n$, $x \in M$ and where $h \in H$ such that
$gg_i=g_jh$.

Now return to our situation. Denote by $L$ the subgroup $\{e,
x^2\}$. Note that there are only two indecomposable modules over
$kL$: the trivial module $k$ and $kL$. If we write $T_0=L \coprod
yL$, then the induced module of $k$ is  $$\text{Ind}_L^{T_0}k=kT_0
\otimes_{kL}k=k(e\otimes 1) \oplus k(y \otimes 1).$$ We obtain
easily
$$\alpha_0 (e \otimes 1) =e \otimes 1+y \otimes 1, \ \ \alpha_0 (y \otimes 1) =y \otimes 1+e \otimes 1$$
$$\beta_0 (e \otimes 1)=e \otimes 1+y \otimes   1,\ \ \beta_0 (y \otimes 1) =y \otimes 1+e \otimes
1$$  If  we write  $f_0=e\otimes 1$ and $f_1=e\otimes 1+y\otimes 1$,
then $\alpha_0 f_0=f_1, \beta_0 f_0=f_1$ and $\alpha_0 f_1=0=\beta_0
f_1$. We obtain the isomorphism
$$\text{Ind}_{L}^{T_0}k=(\xymatrix@1{f_0 \ar@<1ex>[r]^{\alpha_0}
\ar@<-1ex>[r]_{\beta_0}& f_1}) \simeq M(\alpha_0 \beta_0^{-1},1,1)$$
The induction from $\{e, y\}$ to $T_0$ and from $\{e, x^2\}$ to
$T_0$ can be calculated similarly.
\end{Proof}

\bigskip

To prove the main theorem, we will calculate the induced module for
each indecomposable module over $H$, $T_0$ and $T_1$ respectively.

\subsection{Induction from $H$ to $D_8$}

By Lemma~\ref{cyclic}, all indecomposable $kH$-modules are of the
form $M(\gamma^i)$ with $0 \leq i \leq 3$ and the following lemma
computes their induced modules.

\begin{Lem}
\begin{itemize}
\item[(1)] $\text{Ind}_H^{D_8} M(\gamma^0)=\text{Ind}_H^{D_8}
k=M(\beta\alpha^{-1},1,1)$.

 \item[(2)] $\text{Ind}_H^{D_8}
M(\gamma^1)=M(\beta\alpha\beta^{-1}\alpha^{-1},1,1)$.

\item[(3)] $\text{Ind}_H^{D_8}
M(\gamma^2)=M(\beta\alpha\beta\alpha^{-1}\beta^{-1}\alpha^{-1},1,1)$.
\item[(4)] $\text{Ind}_H^{D_8} M(\gamma^3)=\text{Ind}_H^{D_8} kH=kD_8$.
\end{itemize}
\end{Lem}

 The proof uses to the same argument in Proposition~\ref{klein1} and so it is left to the reader.

\subsection{Induction from $T_0$ to $D_8$}

The element $x\in  D_8$ acts via conjugation over $T_0$ et therefore
induces
 an automorphism of $kT_0$, say $\sigma$. We have
$$\sigma(\alpha_0)=x \alpha_0 x^{-1}=1+xyx^{-1}=1+yx^2=\beta_0$$ and

$$\sigma(\beta_0)=x \beta_0 x^{-1}=1+xyx^2x^{-1}=1+y=\alpha_0$$
So the action of $x$ exchanges   $\alpha_0$ and $\beta_0$. Let $M$
be a $kT_0$-module. Denote by $M^\sigma$ the new $kT_0$-module
obtained via $\sigma$. The following lemma deduces immediately from
the argument above
 et from the constructions of string modules and band modules.

\begin{Lem}\label{sigma}
\begin{itemize}
\item[(1)] Let  $C$ be a string and denote by    $C^\sigma$ the new
string by exchanging    $\alpha_0$ and $\beta_0$. Then
$M(C)^\sigma\cong M(C^\sigma)$.

\item[(2)]Let  $C$ be a band.  Then  for all  $n \in
\mathbb{N}$ and  $\lambda \in k^*$, we have  $M(C, n,
\lambda)^\sigma\cong M(C^\sigma, n, \lambda)$.
\end{itemize}
\end{Lem}

As a consequence, for the Auslander-Reiten quiver of $kT_0$, we
obtain
\begin{Prop}
\begin{itemize}
\item[(1)] Each module $M$ in the component of
$\mathbb{Z}\tilde{A}_{12}$ type is  stable by $\sigma$ (i.e.
$M^\sigma \cong M$).

\item[(2)]  $\sigma$ is an isomorphism from   the  homogenous tube $Q(\alpha_0)$
 to the  homogenous tube  $Q(\beta_0)$.

\item[(3)] Given $\lambda \in k^*$, the component
$Q(\alpha_0\beta_0^{-1}, \lambda)$ is  stable by  $\sigma $ if and
only if  $\lambda=1$

\end{itemize}

\end{Prop}

\begin{Proof}
(1) The modules in the component of $\mathbb{Z}\tilde{A}_{12}$ type
are of the form  $\Omega^n(k)=M((\alpha_0\beta_0^{-1})^n)$ or
$\Omega^{-n}(k)=M(\alpha_0^{-1}\beta_0)^n)$ for all  $n\in
\mathbb{N}_0$. These strings   $C$ in this component verify
$C^\sigma=C^{-1}$ and recall that $M(C^{-1}) \cong M(C)$,
Lemma~\ref{sigma} (1) implies the desired result.

(2) As all module in the component  $Q(\alpha_0)$ (resp.
$Q(\beta_0)$) are of the  forme $M(\alpha_0(\beta_0\alpha_0)^n)$
with $n\in \mathbb{N}_0$ (resp. $M(\beta_0(\alpha_0\beta_0)^n)$ with
$n\in \mathbb{N}_0$), then  by the preceding lemma,
$M(\alpha_0(\beta_0\alpha_0)^n)^\sigma\cong
M(\beta_0(\alpha_0\beta_0)^n)$. $\sigma$ thus establishes an
isomorphism between  $Q(\alpha_0)$ and $Q(\beta_0)$.

(3) $M(\alpha_0\beta_0^{-1}, n, \lambda)^\sigma\cong
M(\beta_0\alpha_0^{-1}, n, \lambda) \cong M(\alpha_0\beta_0^{-1}, n,
1/\lambda)$.

\end{Proof}

\begin{Lem}\label{omega}\begin{itemize}
\item[(1)]
$\text{Ind}_{T_0}^{D_8} \Omega^n(k)\cong \Omega^nM(\beta)$, for all
$n \in \mathbb{Z}$.
\item[(2)]
$$\text{Ind}_{T_0}^{D_8}M(\alpha_0)\cong M(\beta\alpha\beta)\cong \text{Ind}_{T_0}^{D_8}M(\beta_0)$$
As a consequence, the homogeneous   tubes $Q(\alpha_0)$ and
$Q(\beta_0)$ are transformed by induction   onto the same
homogeneous
 tube $Q(\beta\alpha\beta)$.

 \end{itemize}

\end{Lem}
\begin{Proof}
Since the indecomposablity theorem of Green(\cite[Theorem
3.13.3]{Ben}) implies that for any $n\in \mathbb{Z}$,
$$\Omega^n(\text{Ind}_{T_0}^{D_8} k)\cong
\text{Ind}_{T_0}^{D_8}(\Omega^nk),$$ it suffices to prove that
$\text{Ind}_{T_0}^{D_8}k\cong M(\beta)$ which is an easy
calculation.

The isomorphism $\text{Ind}_{T_0}^{D_8}M(\alpha_0)\cong
M(\beta\alpha\beta)$ is also simple to prove and is left to the
reader. Recall that for a group $G$ and $H$ a normal subgroup of
$G$, the inertia group of a component of the Auslander-Reiten quiver
of $kH$ is by definition the set of elements of $G$ whose induced
inner automorphisms of $kG$ map this component to itself. As
$\sigma$ transforms $Q(\alpha_0)$ into $Q(\beta_0)$, the inertia
group  of $Q(\alpha_0)$ is $T_0$ and a theorem of S. Kawata
(\cite{K2})
 implies that   induction from $T_0$ to $G$ induces an isomorphism from $Q(\alpha_0)$
(also from  $Q(\beta_0)$) to $Q(\beta\alpha\beta)$

\end{Proof}

\begin{Lem} If $\lambda \in    k-\{0,1\}$,
$\text{Ind}_{T_0}^{D_8}M(\alpha_0\beta_0^{-1},n, \lambda)\cong
M(\beta\alpha\beta\alpha^{-1}, n, \mu)$ with
$\mu=\frac{\lambda}{\lambda^2+1}$. Consequently, the component
$Q(\alpha_0\beta_0^{-1}, \lambda) $ with  $\lambda \in    k-\{0,1\}$
  becomes after  induction the component
$Q(\beta\alpha\beta\alpha^{-1}, \mu)$.
\end{Lem}

\begin{Proof}Write $$M(\alpha_0\beta_0^{-1}, 1, \lambda)\cong
\left(\xymatrix@1{e_0
\ar@<1ex>[r]^{\alpha_0=1}\ar@<-1ex>[r]_{\beta_0=\lambda}&
e_1}\right).$$
 Then $$\alpha_0 e_0=e_1, \alpha_0e_1=0,
 \beta_0e_0=\lambda e_1, \beta_0e_1=0.$$
 Its induced module is
$$\text{Ind}_{T_0}^{D_8} M(\alpha_0\beta_0^{-1}, 1, \lambda)=k(e\otimes
e_0) \oplus k(e\otimes e_1)\oplus k (x\otimes e_0) \oplus k(x\otimes
e_1).$$
 Direct  calculations yield that
$$\alpha (e \otimes e_0) =e\otimes e_1,
 \alpha (e \otimes e_1)=0, \alpha (x\otimes e_0) =\lambda
x\otimes e_1, \alpha (x\otimes e_1) =0, $$ $$\beta (e\otimes e_0)=
e\otimes e_0+ x\otimes e_0+ \lambda x\otimes e_1, \beta (e\otimes
e_1) =e\otimes e_1+ x\otimes e_1, $$
$$\beta (x\otimes e_0) =
x\otimes e_0+ e\otimes e_0+ \lambda e\otimes e_1, \beta (x\otimes
e_1) =x\otimes e_1+ e\otimes e_1.$$

If we impose $$e_0'=e\otimes e_0+\frac{1}{\lambda}x\otimes e_0,\ \
 e_1'=\frac{\lambda+1}{\lambda}(e\otimes e_0+x\otimes e_0)+e\otimes
e_1+\lambda x\otimes e_1, $$
$$e_2'=\frac{\lambda+1}{\lambda}(e\otimes e_1+\ \lambda x\otimes
e_1,\ \ e_3'=\frac{\lambda^2+1}{\lambda}(e\otimes e_1+ x\otimes
e_1),$$ then we can verify that  $\beta e_0'=e_1',\ \alpha
e_1'=e_2',\ \beta e_2'=e_3'\ \text{and}\ \alpha e_0'= \mu e_3'.$
  The statement follows from the diagram:
$$\xymatrix@R1.0cm@C1.0cm{e_0' \ar[r]^-{\beta=1} \ar[d]_{\alpha=\mu}& e_1'
\ar[d]^-{\alpha=1}\\
 e_3'  & e_2'\ar[l]_-{\beta=1}
 } $$

Since $\lambda\neq 1$, $\sigma$ doesn't stabilize the component
$Q(\alpha_0\beta_0^{-1}, \lambda)$ and the inertia group  de
$Q(\alpha_0\beta_0^{-1}, \lambda)$ is $T_0$, the theorem of S.
Kawata cited above implies that the induction from $T_0$ to $D_8$
induces an isomorphism between $Q(\alpha_0\beta_0^{-1}, \lambda)$
and $Q(\beta\alpha\beta\alpha^{-1}, \mu)$.

\end{Proof}

\begin{Prop}\label{longcorde} For any $n\in \mathbb{N}_1$,
$$\text{Ind}_{T_0}^{D_8}M(\alpha_0\beta_0^{-1}, n, 1)=M(C_n, 1, 1)$$
\end{Prop}

The proof of this proposition which is rather complicated is
postponed to the final section
\subsection{Induction from $T_1$ to $D_8$} All the statements in this subsection can be
proved using the same method as in the previous subsection, so we
omit them.

\begin{Lem}\begin{itemize}
\item[(1)]
$\text{Ind}_{T_1}^{D_8} \Omega^n(k)\cong \Omega^nM(\alpha)$, for all
$n \in \mathbb{Z}$.
\item[(2)]
$$\text{Ind}_{T_1}^{D_8}M(\alpha_1)\cong M(\alpha\beta\alpha )\cong \text{Ind}_{T_1}^{D_8}M(\beta_1)$$
As a consequence, the homogeneous   tubes $Q(\alpha_1)$ and
$Q(\beta_1)$ are transformed by induction   onto the same
homogeneous
 tube $Q(\alpha\beta\alpha )$.

 \end{itemize}

\end{Lem}

\begin{Lem} If $\lambda \in    k-\{0,1\}$,
$\text{Ind}_{T_1}^{D_8}M(\alpha_1\beta_1^{-1},n, \lambda)\cong
M(\alpha\beta\alpha\beta^{-1}, n, \mu)$ with
$\mu=\frac{\lambda}{\lambda^2+1}$. Consequently, the component
$Q(\alpha_1\beta_1^{-1}, \lambda) $ with  $\lambda \in    k-\{0,1\}$
  becomes after  induction the component
$Q(\alpha\beta\alpha\beta^{-1}, \mu)$.
\end{Lem}

\begin{Prop}\label{longcorde1}For arbitrary $n\in \mathbb{N}_1$,
$$\text{Ind}_{T_1}^{D_8}M(\alpha_1\beta_1^{-1}, n, 1)=M(D_n, 1, 1)$$
\end{Prop}

\subsection{Proof of the main theorem}

Since we have calculated all the induced modules, we can deduce the
main theorem from the  calculations in the subsections 3.2-3.4,
taking into account the results recalled at the end of the
section~\ref{mr}.

\section{Induction of string modules}\label{formula}

In this section, let $G$ be the dihedral group of order $2^n$ with
$n\geq 3$. Let $M(C)$ be a string over $kT_0$.

A string  $C=\alpha_1 \cdots \alpha_s$ of strictly positive length
is direct (resp. inverse) if all the $
 \alpha_i$ are direct arrows (resp. formal inverses).
Let  $C=C_1C_2\cdots C_n$  where the substrings $C_1, \cdots, C_n$
are direct or inverse  and such that for each  $1\leq i \leq n-1$,
$C_i $ is  direct  (resp. inverse) if and only if   $C_{i+1}$ is
inverse (resp. direct). These substrings $C_i$ are called
\textit{segments} of $C$.

 Let $ C=C_1C_2\cdots C_n$ be a string over $kT_0$ where
the substrings $C_i$ are  its segments.   We use the convention that
$|C_i|=-1$ for $i\leq 0$ or $i\geq n+1$. Let $1\leq i\leq n-1$ and
define a function $\theta: \mathbb{N}_1 \rightarrow \{+1, 0, -1\}$
  as follows:
\[ \theta(s)=\left\{ \begin{array}  {r@{\quad , \quad }l}
                   1&
                   (\text{if \ $ |C_{i-s+1}|>|C_{i+s}|$  and $s$ is odd) or
(if  $ |C_{i-s+1}|<|C_{i+s}|$  and $s$ is even })\\
                   0   &
                   \text{if \ $ |C_{i-s+1}|=|C_{i+s}|$}   \\
                   -1   &
                   \text{(if \ $ |C_{i-s+1}|>|C_{i+s}|$  and $s$ is even) or (if  $ |C_{i-s+1}|<|C_{i+s}|$  and $s$ is
odd)}
\end{array}
            \right.\]
Let $t$ be the first natural number such that $\theta(t)\neq 0$. If
$\theta(t)=1$, then we define $C_i>C_{i+1}$ and if $\theta(t)=-1$,
$C_i<C_{i+1}$.

With this order at hand, we construct a new string over $kG$, say
 $\varphi(C)=
\tilde{C}_1\tilde{C}_2\cdots \tilde{C}_n$ where for all $1\leq i\leq
n$ the $\tilde{C}_i$
 are the segments of $\varphi(C)$ such that
\begin{itemize}
\item[(1)] For all $1\leq i\leq n$,  $\tilde{C}_i\  \text{is
direct (resp. inverse)}  \Longleftrightarrow
   C_i \  \text{is direct (resp. inverse)}$

\item[(2)]
$   |\tilde{C}_i|=\left\{ \begin{array} {r@{ ,  \quad }l}
                   2|C_i|-1   &
                   \text{if} \  C_i<C_{i-1},C_{i+1} \\
                  2|C_i|+1   &
                   \text{if}\      C_i>C_{i-1},C_{i+1}\\
                  2|C_i| & \text{otherwise}
\end{array}
            \right.
$
\item[(3)] Choose  one $i$ such that $1\leq i\leq n$ and $C_i>C_{i-1},C_{i+1}$,  then we impose that  $\tilde{C}_i$ begins
with $\beta$ or $\beta^{-1}$.
\end{itemize}

We can also    construct similarly a new string $\psi(D)$ over $kG$
from a string $D=D_1\cdots D_m$ over $kT_1$. The difference with the
case $kT_0$ is that the last condition becomes

\begin{center}\textit{(3')Choose one  $i$ such that $1\leq i\leq n$ and $D_i>D_{i-1}, D_{i+1}$, then we impose that
  $\tilde{D}_i$ begins
with $\beta$ or $\beta^{-1}$.}\end{center}

\begin{Rem}
\begin{itemize}
\item[(1)] It is easy to see that if for   $1\leq j\leq n$, $C_j>C_{j-1}, C_{j+1}$, then
$\tilde{C}_j$ begins by $\beta$ or $\beta^{-1}$ and ends by $\beta$
or $\beta^{-1}$; if for $1\leq j\leq n$, $C_j<C_{j-1}, C_{j+1}$,
then $\tilde{C}_j$ begins by $\alpha$ or $\alpha^{-1}$ and ends by
$\alpha$ or $\alpha^{-1}$. In particular, the construction of
$\varphi(C)$ is independent of the choice of $C_i$ in the third
condition.

\item[(2)]As we expect that $\text{Ind}_{T_0}^G M(C)\cong
M(\varphi(C))$, at least the new string has  the 'right' length. In
fact, as always $C_0<C_1$ and $C_n>C_{n+1}$, if there are  $t$
segments $C_i$ such that $C_i<C_{i-1}, C_{i+1}$, then there exist
$t+1$ segments $C_i$ such that $C_i>C_{i-1}, C_{i+1}$. We thus have
$|\varphi(C)|=2|C|+1$.
\end{itemize}
\end{Rem}

\begin{Ex} Let $C=C_1\cdots
C_8=\alpha_0^{-1}\beta_0\alpha_0^{-1}\beta_0\alpha_0^{-1}\beta_0\alpha_0^{-1}\beta_0$.

$$\vcenter{\xymatrix@R1.0cm@C1.0cm{  &\ar@{-}[dl]_{\alpha_0} \ar@{-}[dr]_{\beta_0}&&
\ar@{-}[dl]_{\alpha_0}\ar@{-}[dr]_{\beta_0}
 & &\ar@{-}[dl]_{\alpha_0} \ar@{-}[dr]_{\beta_0}&&\ar@{-}[dl]_{\alpha_0} \ar@{-}[dr]_{\beta_0}& \\
             &&&&&&&&&
                                    }}$$
Then $$C_1>C_2<C_3>C_4>C_5<C_6>C_7<C_8$$ and the new string
$\varphi(C)=\tilde{C}_1\cdots \tilde{C}_8$ is of the form:
$$\vcenter
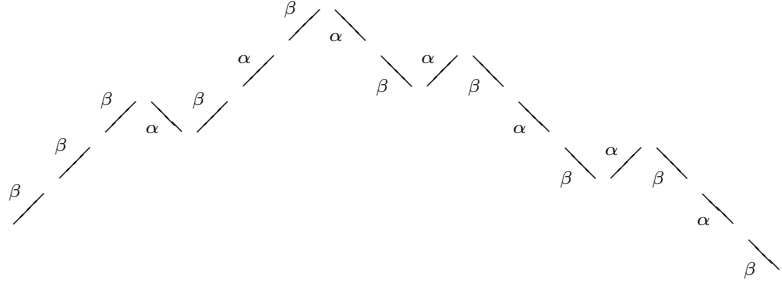
{\xymatrix@R0.4cm@C0.4cm{& & & & & & &
&\ar@{-}[dl]_\beta
\ar@{-}[dr]_\alpha && & & & & & & & & & \\
&& & & & & & \ar@{-}[dl]_\alpha &&\ar@{-}[dr]_\beta
&&\ar@{-}[dl]_\alpha
\ar@{-}[dr]_\beta &&&&&&&\\
& &&&  \ar@{-}[dl]_\beta \ar@{-}[dr]_\alpha
&&\ar@{-}[dl]_\beta &&&&&&\ar@{-}[dr]_\alpha&&&&&&\\
 &&&\ar@{-}[dl]_\beta&&&&&&&&&&\ar@{-}[dr]_\beta&&\ar@{-}[dr]_\beta
 \ar@{-}[dl]_\alpha &&&& \\
             &&\ar@{-}[dl]_\beta&&&&&&&&&&&&&& \ar@{-}[dr]_\alpha&&&\\
             &&&&&&&&&&&&&&&&&\ar@{-}[dr]_\beta&&&\\
             &&&&&&&&&&&&&&&&&&&
  }}$$

\end{Ex}

We give the following
\begin{Conj}

$$\text{Ind}_{T_0}^G M(C)\cong M(\varphi(C))$$ and
$$\text{Ind}_{T_1}^G M(D)\cong M(\psi(D))$$
\end{Conj}

It is easy to verify this conjecture for the dihedral group of order
$8$.

\begin{Prop} The preceding formula   holds when  $G=D_8$ is  the dihedral group
of order $8$.

\end{Prop}
\begin{Proof} As we have calculated the induced module from $T_0$ to $D_8$ for each
string module over $kT_0$ in Lemma~\ref{omega}, we just need to
verify that this is just the string module defined above using the
order. We only consider $\Omega^n(k)$ with $n\geq 1$, similar for
all other string modules.

It is obvious to see that $\Omega^n(k)=M((\alpha_0\beta_0^{-1})^n)$.
We write $(\alpha_0\beta_0^{-1})^n=C_1C_2\cdots C_{2n}$ with the
$C_i$ being its segments.  We now compare  its segments. The result
can be illustrated as follows:
$$C_1>C_2<\cdots C_n>C_{n+1}\cdots >C_{2n-1}<C_{2n}$$
in which the symbols $>$ and $<$ appear in the alternating way from
$C_1$ to $C_n$ with $C_1>C_2$ and from $C_{2n}$ to $C_{n+1}$ with
$C_{2n}>C_{2n-1}$. We thus obtain  the following description of
$\varphi((\alpha_0\beta_0^{-1})^n)$.
\begin{itemize}
\item[(1)]
$\varphi((\alpha_0\beta_0^{-1}))=\tilde{C}_1\tilde{C}_2=(\beta\alpha\beta)(\beta\alpha)^{-1}$.
For $n\geq 1$, $\varphi((\alpha_0\beta_0^{-1})^{2n+1})$ is obtained
by adding $( \beta\alpha\beta)\alpha^{-1}$ to the left side and the
right side of $\varphi((\alpha_0\beta_0^{-1})^{2n-1})$.

\item[(2)] $\varphi((\alpha_0\beta_0^{-1})^2)=\tilde{C}_1\tilde{C}_2\tilde{C}_3\tilde{C}_4=(
\beta\alpha\beta)(\beta\alpha)^{-1}\alpha(\beta\alpha\beta)^{-1}$.
For $n\geq 2$, $\varphi((\alpha_0\beta_0^{-1})^{2n})$ is obtained by
adding $(\beta\alpha\beta)\alpha^{-1}$ to the left side and the
right side of $\varphi((\alpha_0\beta_0^{-1})^{2n-2})$.

\end{itemize}
Now we can verify without difficulty that $\Omega^n(M(\beta))\cong
M(\varphi((\alpha_0\beta_0^{-1})^n))$.

\end{Proof}

To conclude this section, one notices

\begin{Rem} \begin{itemize}
\item[(1)]In \cite{Zhou}, using a constructive method we verified this
formula in the case that there exists at most one segment $C_i$ such
that $C_{i-1}>C_i<C_{i+1}$. The general case remains unsolved.

\item[(2)] If   in general it was true,  iterations of this formula
should give the  vertices of all string modules.

\item[(3)] It will be nice  if we can    extend this formula to band modules.
\end{itemize}
\end{Rem}

\section{Proof of Proposition~\ref{longcorde}}

Before giving the proof of Proposition~\ref{longcorde},   consider
in detail the structure of the bands $C_n$.

 For $n=2$, $$C_2=
\vcenter
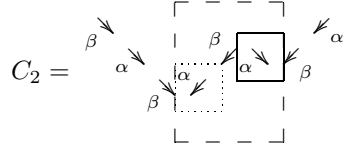
{\xymatrix@R0.2cm@C0.2cm{\ar[dr]_\beta&&&& & & &&
\ar[dl]^\alpha\\
& \ar[dr]_\alpha & & & & \ar[dr]_\alpha \ar[dl]_\beta&& \ar[dl]^\beta& \\
& &  \ar[dr]_\beta& &\ar[dl]_\alpha&&& & & \\
&
 && & & & & & &\\
 &&&&&&&&&
 \save"1,4"."5,7"*[F--]\frm{}\restore
\save"3,4"."4,5"*[F.]\frm{}\restore
\save"2,6"."3,7"*[F-]\frm{}\restore}}$$ Denote by $C_2'$ the part
with  boundary $---$ which is $\alpha^{-1}\beta^{-1}\alpha$, by
$C_2^{(1)}$ the part with boundary $\xymatrix@1{\ar@{.}[r]&}$ which
is $\alpha^{-1} $  and by $C_2^{(2)}$
 the part with boundary $\xymatrix@1{\ar@{-}[r]&}$ which is $ \alpha$.
We see easily that  $(C_2^{(1)})^{-1}$ is equal to
    $C_2^{(2)}$ as strings (in fact $\alpha$). We then have  $C_2=\beta\alpha\beta C_2^{(1)}\beta^{-1}C_2^{(2)}\beta^{-1}
    \alpha^{-1}$.

 For $n=3$,
$$C_3=
\vcenter
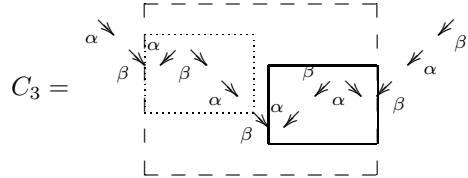
{\xymatrix@R0.2cm@C0.2cm{ \ar[dr]_\alpha& & & & & & & &
& & & & \ar[dl]^\beta\\
& \ar[dr]_\beta & & \ar[dl]_\alpha \ar[dr]_\beta&&&& & & &&
\ar[dl]^\alpha&\\
& & & & \ar[dr]_\alpha & & & & \ar[dr]_\alpha \ar[dl]_\beta&& \ar[dl]^\beta& & \\
& & & & &  \ar[dr]_\beta& &\ar[dl]_\alpha&&& & & & \\&& & &
 && & & & & & & &\\
 &&&&&&&&&&&&&
 \save"1,3"."6,10"*[F--]\frm{}\restore
\save"2,3"."4,6"*[F.]\frm{}\restore
\save"3,7"."5,10"*[F-]\frm{}\restore}}$$ Denote by $C_3'$ the part
with boundary $---$,  by $C_3^{(1)}$ the part with boundary
$\xymatrix@1{\ar@{.}[r]&}$ and by $C_3^{(2)}$
 the part with boundary $\xymatrix@1{\ar@{-}[r]&}$. We see easily that  $(C_3^{(2)})^{-1}$
is equal to
    $C_3^{(1)}$ as strings (in fact, $\alpha^{-1}\beta\alpha $). We
   then  have  $C_3=\alpha\beta C_3^{(1)}\beta
    C_3^{(2)}\beta^{-1}\alpha^{-1}\beta^{-1}$.

If $n$ is even   and  $n\geq 4$, we have $C_n=\beta\alpha\beta
C_n^{(1)}\beta^{-1}C_n^{(2)}\beta^{-1}
    \alpha^{-1}$ with $(C_n^{(1)})^{-1}=C_n^{(2)}$. In fact, by the construction of $C_n$, we can write
    $C_n=\beta\alpha\beta
E_nC_2F_n\beta^{-1}
    \alpha^{-1}$ for certain strings $E_n$ and $F_n$, then we impose $C_n^{(1)}=E_n\beta\alpha\beta\alpha^{-1}$
     and $C_n^{(2)}=\alpha\beta^{-1}\alpha^{-1}F_n$. It is easy  to see by induction that  $(C_n^{(1)})^{-1}=C_n^{(2)}$.
      The situation can be illustrated by
$$
 
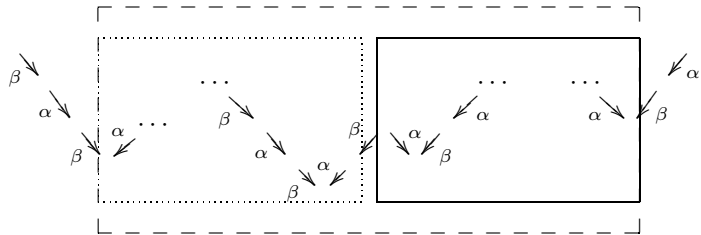
\xymatrix@R0.2cm@C0.2cm{&&&&&&&&&&&&&&&&&&& &&&\\
  \ar[dr]_\beta& & & & & & & & & & & & & &  &   & & & \ar[dl]^\alpha\\
 & \ar[dr]_\alpha& & & & \cdots\ar[dr]_\beta& & & & & & & &
\ar[dl]^\alpha \cdots & &  \cdots\ar[dr]_\alpha & &\ar[dl]^\beta & \\
  & &\ar[dr]_\beta & &\ar[dl]_\alpha \cdots & & \ar[dr]_\alpha& & &
&\ar[dl]_\beta \ar[dr]^\alpha& & \ar[dl]^\beta& & &   & & &  \\
 & & & & & & &\ar[dr]_\beta & & \ar[dl]_\alpha& & & & & &    & & &  \\
 & & & & & & & & & & & & & &  & & & &    \\
 & & & & & & & & & & & & & &  & & & &
\save"1,4"."7,17"*[F--]\frm{}\restore
\save"2,4"."6,10"*[F.]\frm{}\restore
\save"2,11"."6,17"*[F-]\frm{}\restore}$$ where $C_n'$ is the part
with boundary $---$,    $C_n^{(1)}$ is the part with boundary
$\xymatrix@1{\ar@{.}[r]&}$ and   $C_n^{(2)}$ is
 the part with boundary $\xymatrix@1{\ar@{-}[r]&}$. Notice that the string
$C_2$ appears in the middle of   this diagram (and also  in the
middle of  all the diagrams which appear  from now on and which
contain $C_2$).

If $n$ is odd   and  $n\geq 5$, $C_n=\alpha\beta C_3^{(1)}\beta
    C_3^{(2)}\beta^{-1}\alpha^{-1}\beta^{-1}$ with $(C_n^{(2)})^{-1}=C_n^{(1)}$. This can be proved as above.
    The situation can be illustrated  by
$$
 
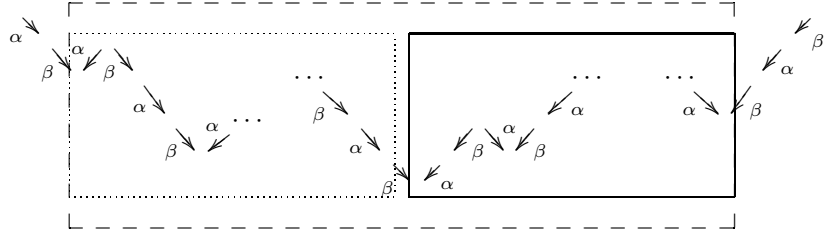
\xymatrix@R0.2cm@C0.2cm{\ar[dr]_\alpha&&&&&&&&&&&&&&&&&&& &&&\ar[dl]^\beta\\
 &\ar[dr]_\beta&&\ar[dl]_\alpha \ar[dr]_\beta& & & & & & & & & & & & & & &   & & & \ar[dl]^\alpha\\
 &&&& \ar[dr]_\alpha& & & & \cdots\ar[dr]_\beta& & & & & & & &\ar[dl]^\alpha \cdots& &
  \cdots\ar[dr]_\alpha & &\ar[dl]^\beta & \\
 &&& & &\ar[dr]_\beta & &\ar[dl]_\alpha \cdots & & \ar[dr]_\alpha& & & &\ar[dl]^\beta
 \ar[dr]^\alpha& & \ar[dl]^\beta& & &   & & &  \\
 &&&& & & & & & &\ar[dr]_\beta & & \ar[dl]^\alpha& & & & & &    & & &  \\
 &&&& & & & & & & & & & & & & &  & & & &    \\
 &&&& & & & & & & & & & & & & &  & & & &
\save"1,3"."7,20"*[F--]\frm{}\restore
\save"2,3"."6,11"*[F.]\frm{}\restore
\save"2,12"."6,20"*[F-]\frm{}\restore}$$ where $C_n'$ is the part
with boundary $---$,    $C_n^{(1)}$ is the part with boundary
$\xymatrix@1{\ar@{.}[r]&}$ and   $C_n^{(2)}$ is
 the part with boundary $\xymatrix@1{\ar@{-}[r]&}$.

\bigskip

We now begin the proof of Proposition~\ref{longcorde}.

Given  $$M(\alpha_0\beta_0^{-1}, n, 1)=( \xymatrix@1{e_1
\ar@<1ex>[r]^{\alpha_0=Id} \ar@<-1ex>[r]_{\beta_0=J_n(1)}& e_2})$$
where $e_1=(e_{11},   \cdots, e_{1n})^t$ and $ e_2=(e_{21}, \cdots,
e_{2n})^t$ and where $Id$ is the  identity matrix of size $n\times
n$ and where $J_n(1)$ is the Jordan block
$$J_n(1)=\left(\begin{array}{cccc} 1& 1& &0   \\& 1&\ddots  &
\\&  & \ddots &1\\ 0 & &  & 1\end{array}\right) $$
We have for all  $1\leq i \leq n$, $\alpha_0 e_{1i}=e_{2i}$,
$\alpha_0 e_{2i}=0$,
 $\beta_0
e_{1i}=e_{2i}+e_{2,i-1}$  and  $\beta_0 e_{2i}=0$  where  we use the
convention that  $e_{2,0}=0$. The induced module $\text{Ind}_{T_0}^G
M(\alpha_0\beta_0^{-1}, n, 1)$ is
$$ (\bigoplus_{i=1}^{n}k e\otimes e_{1i}) \oplus (\bigoplus_{i=1}^nk
x\otimes e_{1i})\oplus (\bigoplus_{i=1}^nk e\otimes e_{2i}) \oplus
(\bigoplus_{i=1}^nk x\otimes  e_{2i})$$ Direct calculations give
that for all  $1\leq i \leq n$, $\alpha (e\otimes e_{1i}) =e\otimes
e_{2i}$,  $ \alpha (x\otimes e_{1i}) =x\otimes e_{2i}+x\otimes
e_{2,i-1}$, $ \alpha (e\otimes e_{2i})=0$, $ \alpha (x\otimes
e_{2i}) =0$, $ \beta (e\otimes e_{1i}) =e\otimes e_{1i}+x\otimes
e_{1i}+x\otimes e_{2i}+x\otimes e_{2,i-1}$, $ \beta (x\otimes
e_{1i}) =x\otimes e_{1i}+e\otimes e_{1i}+e\otimes e_{2i}+e\otimes
e_{2,i-1}$, $ \beta (e\otimes e_{2i}) =e\otimes e_{2i}+x\otimes
e_{2i}$ and $ \beta (x\otimes e_{2i}) =x\otimes e_{2i}+e\otimes
e_{2i}$.

Now we construct    an explicit basis of  $\text{Ind}_{T_0}^G
M(\alpha_0\beta_0^{-1}, n, 1)$   for $n \leq 3$ which establishes
the isomorphism
$$\text{Ind}_{T_0}^G
M(\alpha_0\beta_0^{-1}, n, 1)\cong M(C_n, 1, 1)$$

\textit{Case $n=1$.}
$$ 
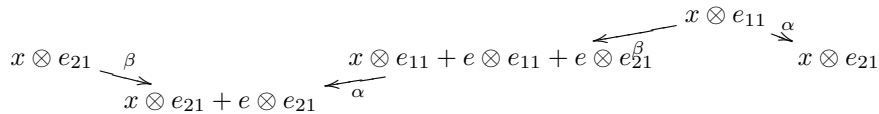
\xymatrix@R0.1cm@C0.2cm{&&& x\otimes e_{11}
\ar[dr]^-\alpha \ar[dl]^-\beta&\\
x\otimes e_{21} \ar[dr]^-\beta& &   x\otimes e_{11}+ e\otimes
e_{11}+e\otimes e_{21}  \ar[dl]^-\alpha&&x\otimes e_{21}\\&
 x\otimes e_{21}+e\otimes e_{21}&& & &} $$
In this diagram, the element in each position is given and it is
easy to see that they  form a basis  of $\text{Ind}_{T_0}^G
M(\alpha_0\beta_0^{-1}, 1, 1)$ (of course, we have to delete  one
$x\otimes e_{21}$). This diagram gives  the desired isomorphism
$\text{Ind}_{T_0}^G M(\alpha_0\beta_0^{-1}, 1, 1)\cong M(C_1, 1,
1)$.

\textit{Case $n=2$.} ( we have   turned the diagram of 90 degrees in
the clockwise direction and we will do this for all diagrams which
appear from now on)
$$
 
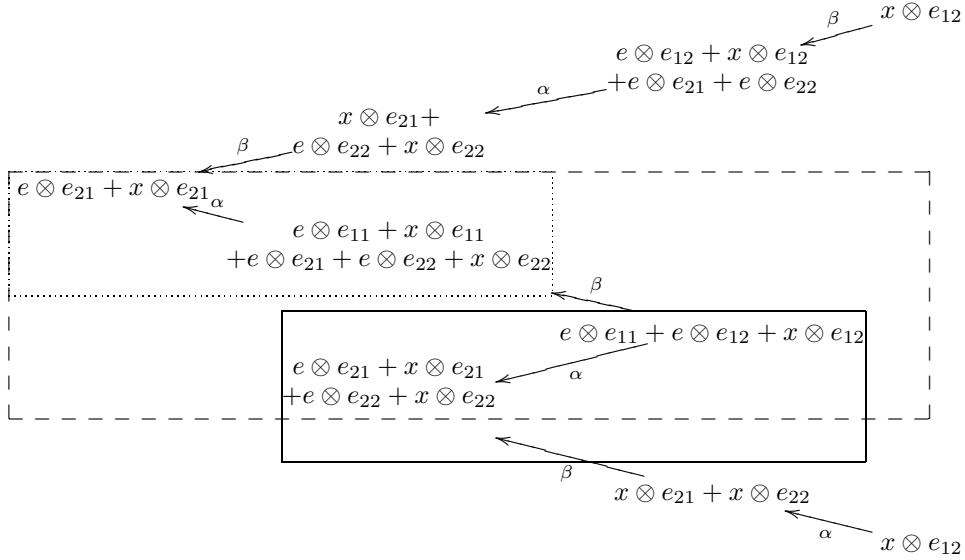
\xymatrix@R0.2cm@C0.1cm{& & &x\otimes
e_{12}\ar[dl]_-\beta \\
& &*\txt{$e\otimes e_{12}+x\otimes e_{12}$\\$+e\otimes
e_{21}+e\otimes
e_{22}$}\ar[dl]_-\alpha& \\
&*\txt{$x\otimes
e_{21}+$\\$e\otimes e_{22}+x\otimes e_{22}$}\ar[dl]_-\beta& &\\
e\otimes e_{21}+x\otimes e_{21}\ar@{<-}[dr]^-\alpha& & &\\
&*\txt{$e\otimes e_{11}+x\otimes e_{11}$\\$+e\otimes e_{21}+e\otimes
e_{22}+x\otimes
e_{22}$\\$   $}\ar@{<-}[dr]^-\beta& & \\
& &*\txt{ $ $\\ $e\otimes e_{11}+e\otimes
e_{12}+x\otimes e_{12}$} \ar[dl]^-\alpha& \\
&*\txt{$e\otimes e_{21}+x\otimes e_{21}$\\$+e\otimes e_{22}
+x\otimes e_{22}$\\$$\\$ $}& & \\
& &x\otimes e_{21}+x\otimes e_{22}\ar@{<-}[dr]_-\alpha\ar[ul]^-\beta& \\
& & &x\otimes e_{12} \save
"4,1"."5,2"*[F.]\frm{}\restore\save"7,4"."4,1"*[F--]\frm{}\restore\save"7,2"."6,3"*[F-]\frm{}\restore}$$
 As
in the case $n=1$, this diagram implies  the isomorphism
$\text{Ind}_{T_0}^G M(\alpha_0\beta_0^{-1}, 2, 1) \cong M(C_2, 1,
1)$. Remark that    the part with boundary  $---$ is  $C_2'$, the
part with boundary    $\xymatrix@1{\ar@{.}[r]&}$ is $C_2^{(1)}$ and
the part with boundary $\xymatrix@1{\ar@{-}[r]&}$ is $C_2^{(2)}$.
Since as  strings, $(C_2^{(1)})^{-1}$ is equal to  $C_2^{(2)}$, if
in the diagram   $C_2'$  we add to the position of  $C_2^{(2)}$ the
diagram
  $(C_2^{(1)})^{-1}$ (with  the elements already given in
$(C_2^{(1)})^{-1}$), then the diagram $C_2' $ becomes  the following
diagram, denoted by $\tilde{C}_2' $,
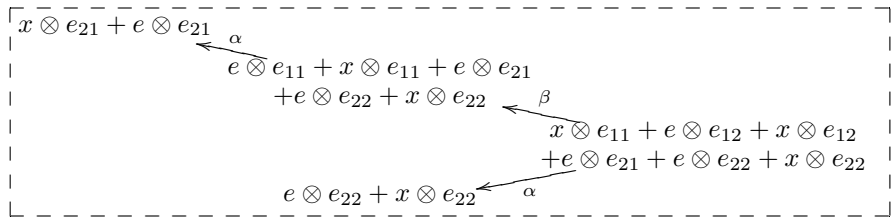
$$
 \xymatrix@R0.2cm@C0.1cm{
x\otimes e_{21}+e\otimes e_{21}\ar@{<-}[dr]^-\alpha& & &&\\
&*\txt{$e\otimes e_{11}+x\otimes e_{11}+e\otimes e_{21}$\\$+e\otimes
e_{22}+x\otimes
e_{22}$}\ar@{<-}[dr]^-\beta& & &\\
& &*\txt{$x\otimes e_{11}+e\otimes
e_{12}+x\otimes e_{12}$\\$ +e\otimes e_{21}+e\otimes e_{22}+x\otimes e_{22}$} \ar[dl]^-\alpha& \\
& *\txt{$e\otimes e_{22} +x\otimes e_{22}$\\$$}  & &&
  \save"1,1"."4,4"*[F--]\frm{}\restore
 }$$

\textit{Case $n=3$} \\
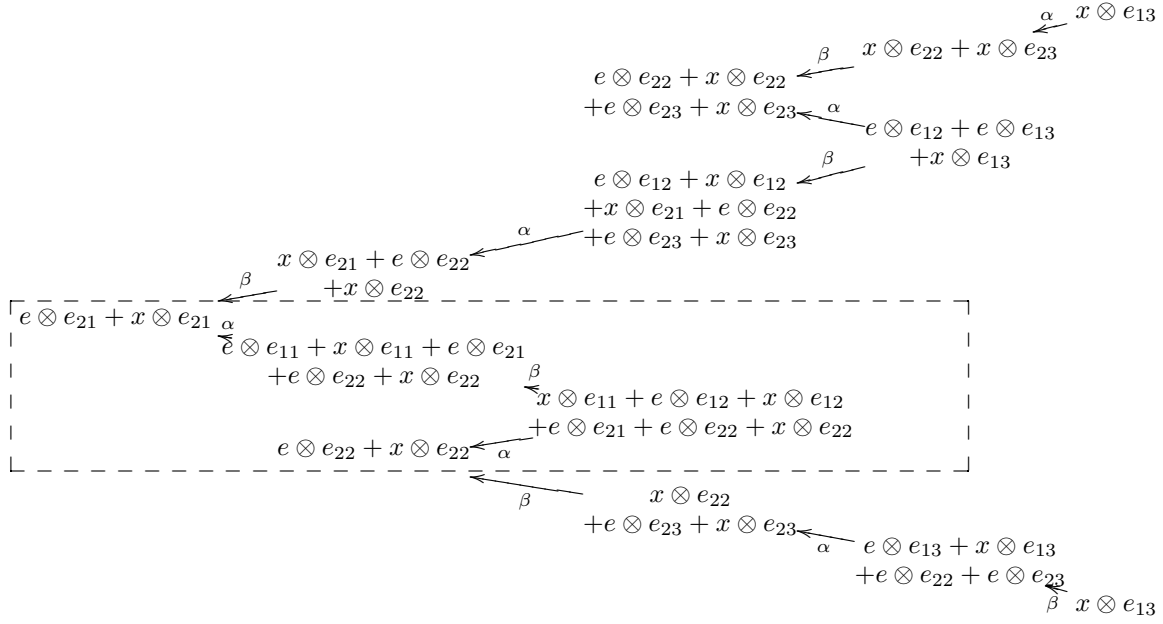
\xymatrix@R0.03cm@C0.01cm{& & & &x\otimes
e_{13}\ar[dl]_-\alpha\\
& & &  x\otimes e_{22}+x\otimes e_{23}\ar[dl]_-\beta & \\
&& *\txt{$e\otimes e_{22}+x\otimes e_{22}$\\$+e\otimes
e_{23}+x\otimes
e_{23}$}\ar@{<-}[dr]^-\alpha& &  \\
& & &*\txt{$e\otimes e_{12}+e\otimes e_{13}$\\$+x\otimes
e_{13}$}\ar[dl]_-\beta& \\
& &*\txt{$e\otimes e_{12}+x\otimes e_{12}$\\$+x\otimes
e_{21}+e\otimes e_{22} $\\$+e\otimes
e_{23}+x\otimes e_{23}$}\ar[dl]_-\alpha & & \\
&*\txt{$x\otimes e_{21}+e\otimes e_{22}$\\$+x\otimes e_{22}$}\ar[dl]_-\beta& & & \\
e\otimes e_{21}+x\otimes e_{21}\ar@{<-}[dr]^-\alpha& & & & \\
&*\txt{$e\otimes e_{11}+x\otimes e_{11}+ e\otimes e_{21}$\\
$+e\otimes e_{22}+x\otimes
e_{22}$}\ar@{<-}[dr]^-\beta& & &\\
&& *\txt{$x\otimes e_{11}+e\otimes e_{12}+x\otimes
e_{12}$\\$+e\otimes e_{21}+e\otimes e_{22}+x\otimes e_{22}$}
\ar[dl]^-\alpha& & \\
& *\txt{$e\otimes e_{22}+x\otimes e_{22}$\\$$}\ar@{<-}[dr]_-\beta& & &\\
&&*\txt{$x\otimes e_{22}$\\$+e\otimes e_{23}+x\otimes e_{23}$}\ar@{<-}[dr]_-\alpha& & \\
&&&*\txt{$e\otimes e_{13}+x\otimes e_{13}$\\$ +e\otimes
e_{22}+e\otimes
e_{23}$}\ar@{<-}[dr]_-\beta& \\
& &&  & x\otimes e_{13}\save "7,1"."10,4"*[F--]\frm{}\restore}\\
We see easily that this diagram gives the desired isomorphism.
Remark that the part in the box, which is equal to  $C_2'$ as
strings, is exactly  the diagram   $\tilde{C}_2'$.

The induction hypothesis for $n-1\geq 3$ is the following:

 (1) $\text{Ind}_{T_0}^GM(\alpha_0\beta_0^{-1}, n-1, 1) \cong
M(C_{n-1}, 1,1)$

(2)There exists a basis of
$\text{Ind}_{T_0}^GM(\alpha_0\beta_0^{-1}, n-1, 1)$ which gives the
isomorphism and which contains the elements already given in the
following  diagrams:

(i) If  $n-1$ is odd and $n-1\geq 3$, the basis of
$\text{Ind}_{T_0}^GM(\alpha_0\beta_0^{-1}, n-1, 1)$ is of the form\\
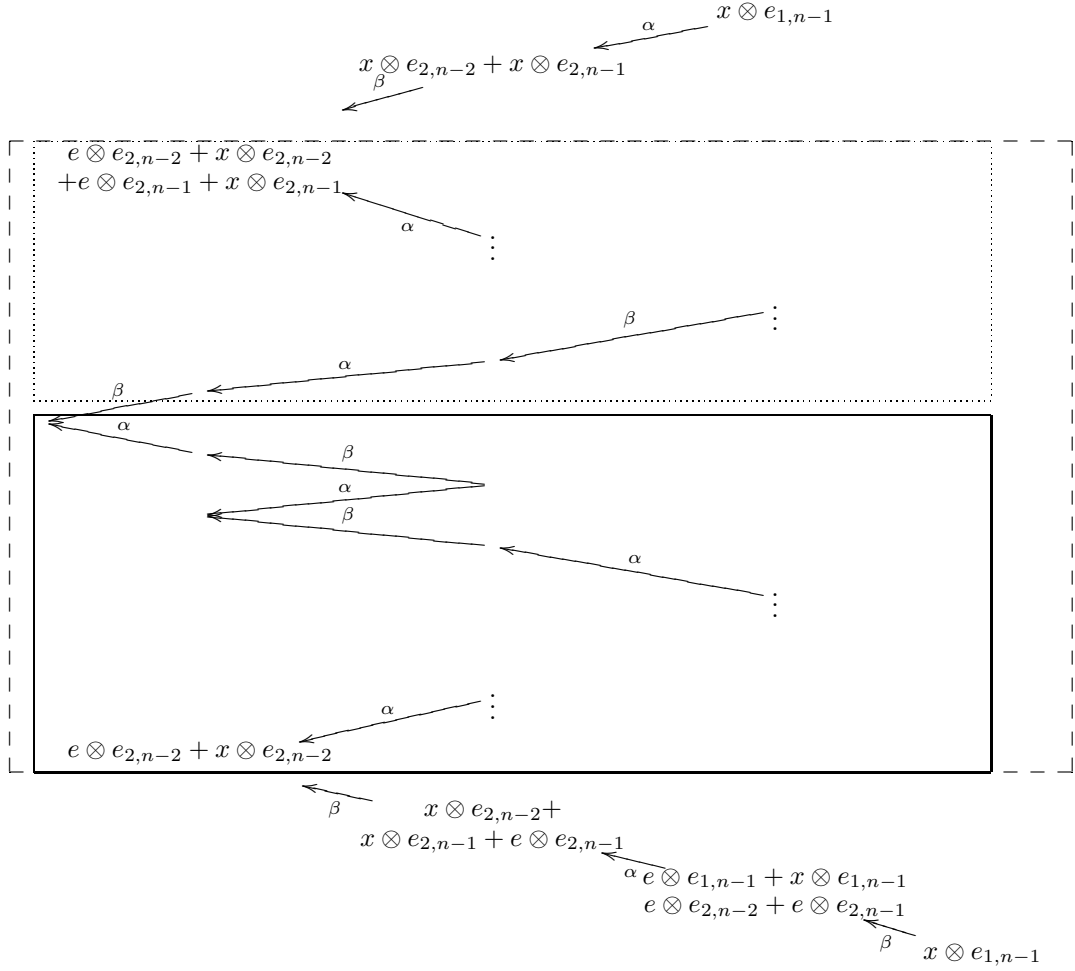
\xymatrix@R0.2cm@C0.1cm{
&&&&x\otimes e_{1,n-1}\ar[dl]_-\alpha&&\\
&&&x\otimes e_{2,n-2}+x\otimes e_{2,n-1}\ar[dl]_-\beta&&&\\
&&*\txt{$ $\\$ $\\$e\otimes e_{2,n-2}+x\otimes e_{2,n-2}$\\ $+e\otimes e_{2,n-1}+x\otimes e_{2,n-1}$}&&&\\
&&&\ar[ul]^-\alpha\vdots&&&\\
  && &&\vdots\ar[dl]_-\beta&&\\
&&&\ar[dl]_-\alpha&&&\\
&&\ar[dl]_-\beta&&&&\\
&&&&&&\\
&&\ar[ul]_-\alpha&&&&\\
&&&\ar[ul]_-\beta\ar[dl]_-\alpha&&&\\
&&&&&&\\
&&&\ar[ul]_-\beta&&&\\
&&&&\ar[ul]_-\alpha \vdots&&\\
&&&&&&\\
&&&\vdots\ar[dl]_-\alpha&&&&\\
&&*\txt{$e\otimes e_{2,n-2}+x\otimes e_{2,n-2}$\\$ $}&&&&&\\
&&&*\txt{$x\otimes e_{2,n-2}+$\\$x\otimes e_{2,n-1}+e\otimes e_{2,n-1}$}\ar[ul]^-\beta&&&&\\
&&&&*\txt{$e\otimes e_{1,n-1}+x\otimes e_{1,n-1}$\\$e\otimes e_{2,n-2}+e\otimes e_{2,n-1}$}\ar[ul]^-\alpha&&&\\
&&&&&x\otimes e_{1,n-1}\ar[ul]^-\beta&&\\
&&&&&&&\\
 \save"3,1"."16,7"*[F--]\frm{}\restore
\save"3,2"."7,6"*[F.]\frm{}\restore
\save"8,2"."16,6"*[F-]\frm{}\restore}

 The part with boundary $\xymatrix@1{\ar@{.}[r]&}$ is
$C_{n-1}^{(1)}$, the part with boundary  $\xymatrix@1{\ar@{-}[r]&}$
is  $C_{n-1}^{(2)}$ and the part with boundary   $---$ is
$C_{n-1}'$. Since as strings, $(C_{n-1}^{(2)})^{-1}$ is equal to
$C_{n-1}^{(1)}$, if we add in $C_{n-1}'$ to the position of
$C_{n-1}^{(1)}$ the diagram $(C_{n-1}^{(2)})^{-1}$ (with the given
elements), then the diagram $C_{n-1}'$ becomes a diagram, denoted by
$\tilde{C}_{n-1}'$, whose 'highest' element is   $e\otimes
e_{2,n-1}+x\otimes e_{2,n-1} $ and whose 'lowest' element is
$e\otimes e_{2,n-2}+x\otimes e_{2,n-2}$.

(ii) If  $n-1$ is even and $n-1\geq 4$, the basis of
$\text{Ind}_{T_0}^GM(\alpha_0\beta_0^{-1}, n-1, 1)$ is of the form
$$
 \xymatrix@R0.6cm@C0.2cm{
 && & & &x\otimes e_{1,n-1}\ar[dl]^-\beta&\\
&&&&*\txt{$e\otimes e_{1,n-1}+x\otimes e_{1,n-1}$\\$+e\otimes
e_{2,n-2}+e\otimes e_{2,n-1}$}
\ar[dl]^-\alpha&&\\
&&&*\txt{$x\otimes e_{2,n-2}+$\\$e\otimes e_{2,n-1}+x\otimes e_{2,n-1}$}\ar[dl]^-\beta&&&\\
&&*\txt{$ $\\$e\otimes e_{2,n-2}+x\otimes e_{2,n-2}$}&&&\\
&&&\ar[ul]^-\alpha\vdots&&&\\
&  & &&\vdots\ar[dl]^-\beta&&\\
&&&\ar[dl]^-\alpha&&&\\
&&\ar[dl]^-\beta&&&&\\
&&&&&&\\
&&\ar[ul]^-\alpha&&&&\\
&&&\ar[ul]^-\beta\ar[dl]^-\alpha&&&\\
&&&&&&\\
&&&\ar[ul]^-\beta&&&\\
&&&&\ar[ul]^-\alpha\vdots&&\\
&&&&& &\\
&&&& \ar[dl]^-\alpha\vdots&&\\
&&&*\txt{$e\otimes e_{2,n-2}+x\otimes e_{2,n-2}$\\$+e\otimes e_{2,n-1}+x\otimes e_{2,n-1}$\\$ $\\$ $}&&&\\
&&&&x\otimes e_{2,n-2}+x\otimes e_{2,n-1}\ar[ul]^-\beta&&\\
&&&&&x\otimes e_{1,n-1}\ar[ul]^-\alpha
\save"4,1"."17,7"*[F--]\frm{}\restore
\save"4,2"."10,6"*[F.]\frm{}\restore
\save"11,2"."17,6"*[F-]\frm{}\restore}$$
 The part with boundary
$\xymatrix@1{\ar@{.}[r]&}$ is $C_{n-1}^{(1)}$, the part with
boundary  $\xymatrix@1{\ar@{-}[r]&}$ is  $C_{n-1}^{(2)}$ and the
part with boundary   $---$ is $C_{n-1}'$. Since as strings,
$(C_{n-1}^{(1)})^{-1}$ is equal to  $C_{n-1}^{(2)}$, if we add in
$C_{n-1}'$ to the position of  $C_{n-1}^{(2)}$ the diagram
$(C_{n-1}^{(1)})^{-1}$ (with the given elements), then the diagram
$C_{n-1}'$ becomes a diagram, denoted by  $\tilde{C}_{n-1}'$, whose
'highest' element is  $e\otimes e_{2,n-2}+x\otimes e_{2,n-2}$ and
whose 'lowest' element is  $e\otimes e_{2,n-1}+x\otimes e_{2,n-1} $.

This finishes the statement of the induction hypothesis.

We now construct the diagram $C_n$ which establishes the desired
isomorphism.

If $n$ is even and $n\geq 4$,    at first we construct  an
incomplete diagram which is $C_n$ as a string and which contains
some given elements.
$$
 
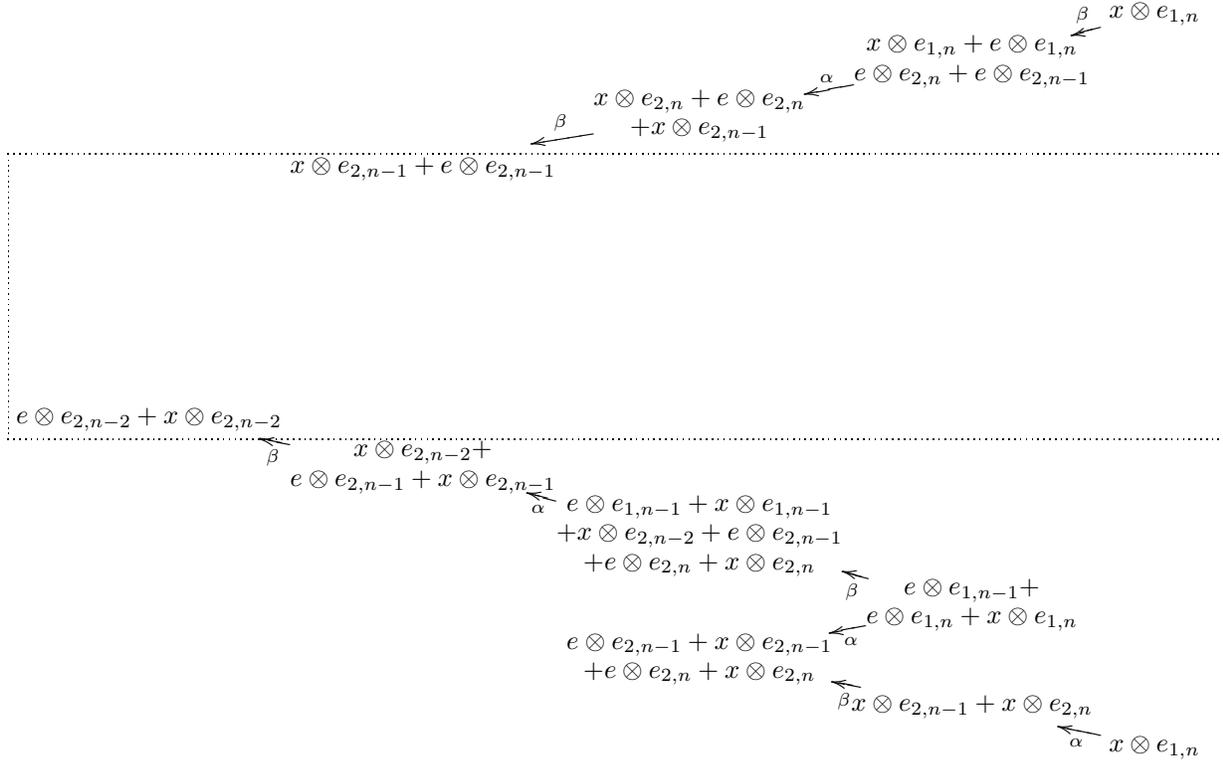
\xymatrix@R0.03cm@C0.01cm{& & & &x\otimes e_{1,n}\ar[dl]_-\beta&\\
&&&*\txt{$x\otimes e_{1,n }+e\otimes e_{1,n }$\\$e\otimes e_{2,n }+e\otimes e_{2,n-1}$}\ar[dl]_-\alpha&&\\
&&*\txt{$x\otimes e_{2,n }+e\otimes e_{2,n }$\\$+x\otimes e_{2,n-1}$}\ar[dl]_-\beta&&&\\
&*\txt{$ $\\$x\otimes e_{2,n-1}+e\otimes e_{2,n-1}$}&&&&&&&\\
&& &&&\\
  & && &&\\
&& &&&\\
& &&&&\\
&&&&&\\
& &&&&\\
&& &&&\\
&&&&&\\
&& &&&\\
&&& &&\\
&&&& &\\
& &&&&\\
e\otimes
e_{2,n-2}+x\otimes e_{2,n-2}&&&&&\\
&*\txt{$x\otimes
e_{2,n-2}+$\\$e\otimes e_{2,n-1}+x\otimes e_{2,n-1}$}\ar[ul]^-\beta&&&&\\
&&*\txt{$e\otimes e_{1,n-1}+x\otimes e_{1,n-1}$\\$+x\otimes
e_{2,n-2}+e\otimes e_{2,n-1} $\\$+e\otimes e_{2,n}+x\otimes e_{2,n}$
}
\ar[ul]^-\alpha&&&\\
&&&*\txt{$e\otimes e_{1,n-1}+$\\$e\otimes e_{1,n}+x\otimes e_{1,n}$}\ar[ul]^-\beta\ar[dl]^-\alpha&&\\
&&*\txt{$e\otimes e_{2,n-1}+x\otimes e_{2,n-1}$\\$+e\otimes e_{2,n }+x\otimes e_{2,n }$}&&&\\
&&&x\otimes e_{2,n-1}+x\otimes e_{2,n }\ar[ul]^-\beta&&\\
&&&&x\otimes e_{1,n}\ar[ul]^-\alpha\\&&&&&
 \save"4,6"."17,1"*[F.]\frm{}\restore}
$$
Since the empty box  is equal to  $C_{n-1}'$ as a string, we replace
the box by  the diagram  $\tilde{C}_{n-1}'$ constructed in the
induction hypothesis   and it is easy to see that $\tilde{C}_{n-1}'$
glues with the elements already given. We verify that the complete
diagram constructed above gives the desired isomorphism
$\text{Ind}_{T_0}^GM(\alpha_0\beta_0^{-1}, n, 1)\cong M(C_n,1,1)$
and thus satisfied the induction hypothesis.

If $n$ is odd and $n\geq 5$, as above we construct   an incomplete
diagram which is $C_n$ as a string and which contains some  given
elements.

 
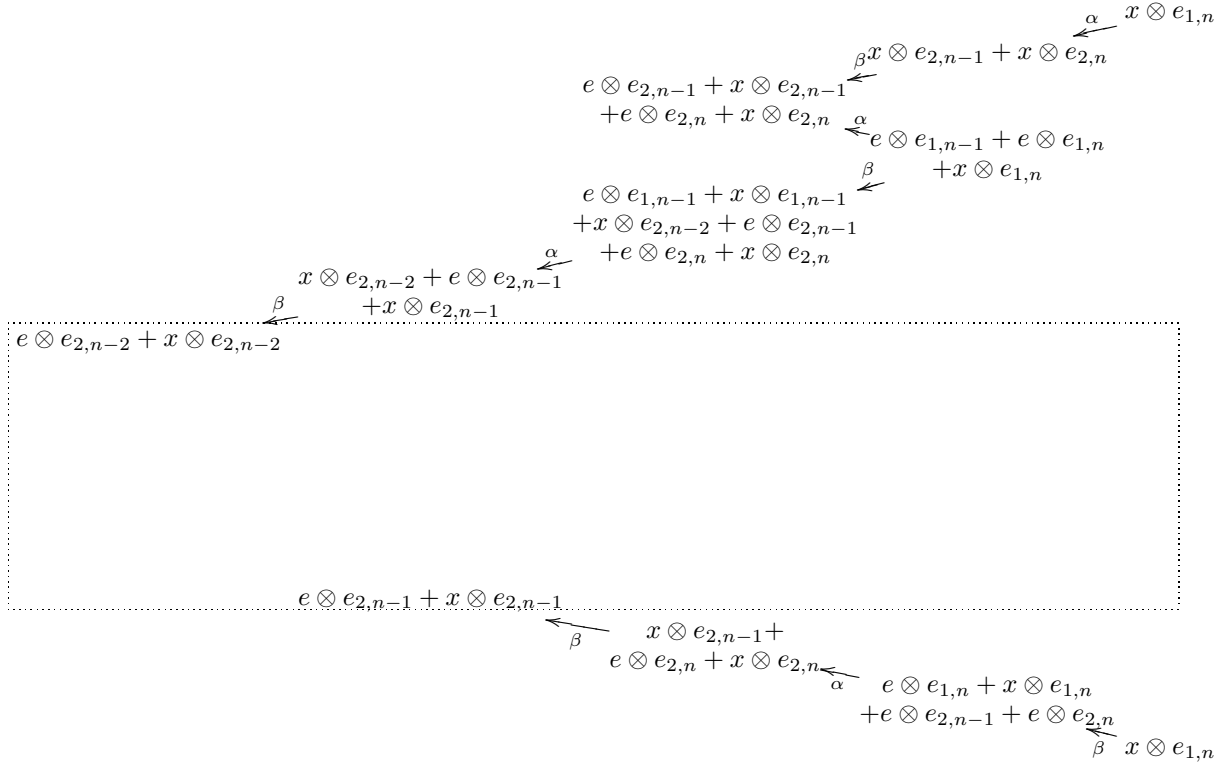
\xymatrix@R0.03cm@C0.01cm{
&&&&x\otimes e_{1,n}\ar[dl]_-\alpha\\
&&&x\otimes e_{2,n-1}+x\otimes e_{2,n }\ar[dl]_-\beta&&\\
&&*\txt{$e\otimes e_{2,n-1}+x\otimes e_{2,n-1}$\\$+e\otimes e_{2,n }+x\otimes e_{2,n }$}&&&\\
&&&*\txt{$e\otimes e_{1,n-1}+e\otimes e_{1,n}$\\$+x\otimes e_{1,n}$}\ar[ul]_-\alpha\ar[dl]_-\beta&&\\
&&*\txt{$e\otimes e_{1,n-1}+x\otimes e_{1,n-1}$\\$+ x\otimes
e_{2,n-2}+e\otimes
e_{2,n-1}$\\$+e\otimes e_{2,n}+x\otimes e_{2,n} $}\ar[dl]_-\alpha&&&\\
&*\txt{$x\otimes
e_{2,n-2}+e\otimes e_{2,n-1}$\\$+x\otimes e_{2,n-1}$}\ar[dl]_-\beta&&&&\\
e\otimes
e_{2,n-2}+x\otimes e_{2,n-2}&&&&&\\
& &&&&\\
  & && &&\\
&& &&&\\
& &&&&\\
&&&&&\\
& &&&&\\
&& &&&\\
&&&&&\\
&& &&&\\
&&& &&\\
&&&&  &\\
&& &&&&\\
&e\otimes e_{2,n-1}+x\otimes e_{2,n-1}&&&\\
&&*\txt{$x\otimes e_{2,n-1}+$\\$e\otimes e_{2,n }+x\otimes e_{2,n }$}\ar[ul]^-\beta&&&\\
&&&*\txt{$e\otimes e_{1,n }+x\otimes e_{1,n}$\\$+e\otimes
e_{2,n-1}+e\otimes e_{2,n }
$}\ar[ul]^-\alpha&&\\
& & & &x\otimes e_{1,n}\ar[ul]^-\beta&\\
&&&&&
 \save"7,1"."20,5"*[F.]\frm{}\restore}
Since the empty box  is equal to  $C_{n-1}'$ as a string, we replace
the box by  the diagram  $\tilde{C}_{n-1}'$ constructed in the
induction hypothesis  and it is easy to see that $\tilde{C}_{n-1}'$
glues with the elements already given. We verify that the complete
diagram constructed above gives the desired isomorphism
$\text{Ind}_{T_0}^GM(\alpha_0\beta_0^{-1}, n, 1)\cong M(C_n,1,1)$
and thus satisfied the induction hypothesis.

This finishes the proof.




\end{document}